\documentclass[reqno,12pt]{amsart}
\usepackage{amsthm}
\usepackage{latexsym,times,paralist}
\usepackage{amssymb,amsmath,amsfonts,esint}
\usepackage{enumitem}
\usepackage[utf8]{inputenc}
\usepackage{color}
	\usepackage{centernot}
\usepackage{hyperref}

\renewcommand{\H}{\mathbb{H}}
\newcommand{\B}{\mathbb{B}}

\newcommand{\G}{\mathbb{G}}

\newcommand{\N}{\mathbb{N}}

\newcommand{\R}{\mathbb{R}}

\newcommand{\T}{\mathbb{T}}
\newcommand{\V}{\mathbb{V}}
\newcommand{\W}{\mathbb{W}}

\renewcommand{\P}{\mathbb{P}}

\newcommand{\cD}{\mathcal{D}}

\newcommand{\cG}{\mathcal{G}}
\newcommand{\cH}{\mathcal{H}}

\newcommand{\cL}{\mathcal{L}}

\newcommand{\cS}{\mathcal{S}}

\newcommand{\cV}{\mathcal{V}}

\newcommand{\cg}{\mathfrak{g}}

\newcommand{\ep}{\varepsilon}

\newcommand{\ph}{\varphi}

\newcommand{\sm}{\setminus}

\newcommand{\lb}{{\big\lbrace}}
\newcommand{\rb}{{\big\rbrace}}

\newcommand{\diam}{\mbox{\rm diam}}

\renewcommand{\exp}{\mbox{\rm exp}\;\!}

\newcommand{\supp}{\mbox{\rm supp}}

\newcommand{\Tan}{\mbox{\rm Tan}}

\newcommand{\Lie}{\mathrm{Lie}}

\newcommand{\q}{\mathrm q}
\newcommand{\m}{\mathrm m}

\newcommand{\h}{\mathrm{h}}

\newcommand{\bcup}{\bigcup}

\newcommand{\res}{\mbox{\LARGE{$\llcorner$}}}

\newcommand{\beqas}{\begin{eqnarray*}}
	\newcommand{\eeqas}{\end{eqnarray*}}
\newcommand{\beqa}{\begin{eqnarray}}
	\newcommand{\eeqa}{\end{eqnarray}}
\newcommand{\beq}{\begin{equation}}
	\newcommand{\eeq}{\end{equation}}
\newcommand{\bce}{\begin{center}}
	\newcommand{\ece}{\end{center}}

\newcommand{\pa}[1]{\left( #1 \right)}               % (  )
\newcommand{\set}[1]{\left\{ #1 \right\}}            % {  }
\newcommand{\qandq}{\quad\mbox{and}\quad}

\newtheorem{The}{Theorem}[section]
\newtheorem{Lem}[The]{Lemma}
\newtheorem{Def}[The]{Definition}
\newtheorem{Rem}[The]{Remark}
\newtheorem{Pro}[The]{Proposition}
\newtheorem{Cor}[The]{Corollary}
\newtheorem{Exa}[The]{Example}
\newtheorem{Con}{Conjecure}

\newcommand{\bt}{\begin{The}}
	\newcommand{\et}{\end{The}}
\newcommand{\bl}{\begin{Lem}}
	\newcommand{\el}{\end{Lem}}
\newcommand{\bd}{\begin{Def}\rm}
	\newcommand{\ed}{\end{Def}}
\newcommand{\br}{\begin{Rem}\rm}
	\newcommand{\er}{\end{Rem}}
\newcommand{\bpr}{\begin{Pro}}
	\newcommand{\epr}{\end{Pro}}
\newcommand{\bc}{\begin{Cor}}
	\newcommand{\ec}{\end{Cor}}
\newcommand{\bj}{\begin{Con}}
	\newcommand{\ej}{\end{Con}}
\newcommand{\bex}{\begin{Exa}}
	\newcommand{\eex}{\end{Exa}}

\theoremstyle{definition}

\numberwithin{equation}{section}

\begin{document}

\title{Characterizations of $k$-rectifiability in homogenous groups}

\author{Kennedy Obinna Idu}
\address{Dip.to di Matematica, Universit\`a di Pisa, Largo Bruno Pontecorvo 5, 56127, Pisa, Italy}
\email{idu@mail.dm.unipi.it}
\author{Valentino Magnani}
\address{Dip.to di Matematica, Universit\`a di Pisa, Largo Bruno Pontecorvo 5, 56127, Pisa, Italy}
\email{valentino.magnani@unipi.it}
\author{Francesco Paolo Maiale}
\address{Scuola Normale Superiore, Piazza dei Cavalieri 7, 56126 Pisa, Italy}
\email{francesco.maiale@sns.it}
\date{\today}
\thanks{V.M. is supported by the University of Pisa, Project PRA 2018 49.}
\date{\today}

\begin{abstract}
A well known notion of $k$-rectifiable set can
be formulated in any metric space using Lipschitz images of subsets
of $\R^k$. We prove some characterizations of $k$-rectifiability,
when the metric space is an arbitrary homogeneous group.
In particular, we show that the a.e.\ existence of the 
$(k,\G)$-approximate tangent group implies $k$-rectifiability.
\end{abstract}

\maketitle

\tableofcontents

\section{Introduction}

Rectifiability is a central concept in Geometric Measure Theory and in many other areas of Mathematics, like Calculus of Variations and Geometric Analysis. It can be naturally defined in metric spaces, considering countably many Lipschitz images of subsets of the Euclidean space $\R^k$, \cite[3.2.14]{Federer69}. We say that the union of these images is a $k$-rectifiable set.

A rich literature is available about characterizations and basic properties of $k$-rectifiable sets in Euclidean spaces and also general metric spaces. For more information, we address the reader to \cite{AmbKir2000Rect,Federer69,Federer1978Colloquium,DeLellis08,mattila}. Actually, more references could be added.

Since the example of Ambrosio and Kirchheim, \cite{AmbKir2000Rect}, it is well known that rectifiability through Lipschitz images does not always work in all metric spaces. The authors showed that the three dimensional Heisenberg group is purely $k$-unrectifiable with $k=2,3,4$. 

A natural notion of ``intrinsic rectifiability'' in codimension one was first proposed in the seminal paper \cite{FSSC01}, where Franchi, Serapioni and Serra Cassano established the De Giorgi rectifiability theorem for finite perimeter sets in all Heisenberg groups.
More generally they introduced $\G$-rectifiable sets as countable unions of level sets of functions with nonvanishing and continuous horizontal gradient in a stratified group $\G$, up to a suitable negligible set, \cite{FSSC5}. 

A number of recent results arose in connection with characterizations of $\G$-rectifiable sets, somehow confirming that they constitute the natural class of ``intrinsic rectifiable sets of codimension one". Merlo showed that in Heisenberg groups such type of rectifiability can be characterized in terms of densities and intrinsically flat tangents, finally reaching the analogous of Preiss' theorem in Heisenberg groups, \cite{Mer20}. 
Although $\G$-rectifiable sets can be very far from being $k$-rectifiable, \cite{KirSer04}, one may still ask whether they can be covered by countable unions of Lipschitz images of more general homogeneous groups, in place of  $\R^k$. These unions of Lipschitz images correspond to the extended notion rectifiability 
introduced by Pauls, \cite{Pau04}. 
Pauls rectifiability and $\G$-rectifiability are in general different, as 
showed by Antonelli and Le Donne in \cite{ALD20}, but their equivalence in Heisenberg groups is still an interesting open question. Recently, the work of Di Donato, F\"assler and Orponen showed that ``intrinsic $C^{1,\alpha}$'' hypersurfaces
can be seen as countable Lipschitz images of 1-codimensional vertical subgroups, \cite{DFO20pr}.

A third notion of rectifiability can be provided through the notion of intrinsic Lipschitz graph, \cite{FSSC06,FSSC11}, considering countable unions of intrinsic Lipschitz graphs. The results of Franchi, Marchi, Serapioni and Serra Cassano essentially showed that this rectifiability corresponds to $\G$-rectifiability as soon as we have a De Giorgi rectifiability theorem, \cite{FMS14,FSSC11}. In Heisenberg groups and for all codimensions, this question has been recently settled by Vittone, using the theory of currents, that lead to a Rademacher type theorem for intrinsic Lipschitz functions, \cite{Vit20pr}. In connection with the De Giorgi's rectifiability problem for sets of finite perimeters, a fourth notion of rectifiability appeared in \cite{DLDMV19pr}.
Still many questions are to be investigated to find a unified view of ``intrinsic rectifiability" in homogeneous groups.

In the present paper, we wish to study all subsets of homogeneous groups that are captured by the classical notion of $k$-rectifiability. 
A first work in this direction is \cite{Pau04}, where the author treats rectifiable sets in an extended sense, considering subsets of stratified groups as domains of the Lipschitz parametrizations. On one hand, a general notion of rectifiability is studied and some first partial results are established. On the other hand, 
the use of Euclidean projections to define cones seems to constitute an obstacle in reaching a characterization of this rectifiability through approximate tangent cones. A stronger notion of approximate tangent cone was needed to prove a rectifiability criterion \cite[Theorem~B]{Pau04}.

The natural notion of ``approximate tangent group" was first introduced in the work of Mattila, Serapioni and Serra Cassano, \cite{MSSC10}, where the projections to define the intrinsic cones arised from the semidirect factorizations of Heisenberg groups.
From the a.e.\ existence of approximate tangent groups in suitable Grassmannians, $k$-rectifiable sets in all Heisenberg groups $\H^n$, with $k\le n$, were characterized, 
along with low codimensional rectifiable sets with positive lower density, \cite{MSSC10}.

Our main result is the characterization of $k$-rectifiable sets in any homogeneous group.

\bt \label{thm.principale}
Let $\G$ be a homogeneous group and let $1\le k \le \upsilon$ be an integer, where $\upsilon$ is introduced in Definition~\ref{d:abelian dimension}. Let $E\subset \G$ be a Borel set such that $\cH^k\res E$ is locally finite.
Then the following conditions are all equivalent:
\begin{enumerate}[label=(\roman*), itemsep=.5em]
	\item The set $E$ is $k$-rectifiable.
	\item For $\cH^k$-a.e. $p\in E$, there exists $\T_p\in \cH(\G, k)$ such that
	\begin{equation} \label{eq.blowup}
	\frac{1}{r^k}(T_{p,r})_\#\cH^k \res E \rightharpoonup \cH^k\res \T_p.
	\end{equation}
	\item For $\cH^k$-a.e. $p\in E$, there exists $\T_p\in \cH(\G, k)$ such that 
	\[\Tan(\cH^k\res E,  p)=\left\{\lambda\cH^k \res \T_p  :  0<\lambda<\infty\right\}.\]
	\item For $\cH^k$-a.e. $p\in E$, there exists $\T_p\in \cH(\G, k)$ such that
	\[ \rm{ap}\Tan_\G^k(E, p) = \T_p. \]
\end{enumerate}
\et
The first important implication is from (i) to (ii), where the crucial role is played
by the ``bi-Lipschitz linearization lemma", stated in Lemma~\ref{lemma.3.1}.
It is a consequence of the a.e.\ differentiability of Lipschitz mappings taking values in homogeneous groups \cite[Theorem~1.1]{MagnaniRajala2014} and the associated area formula, \cite[Theorem~1.2]{MagnaniRajala2014}. These results allow for characterizing $k$-rectifiable sets in the large class of homogeneous groups, where the Lie bracket generating condition on the first layer of their Lie algebra is not required.
The main implication is from (iv) to (i), namely the a.e.\ existence of the $(k,\G)$-approxi\-mate tangent group implies $k$-rectifiability. The proof of this step in \cite{MSSC10} is established by the important positive lower density theorem, see \cite[Theorem 3.10]{MSSC10}.

For the same implication, we follow a different method, fixing our attention on the purely $k$-unrectifiable part of the given set. Our main tool is Theorem~\ref{ppwpod}, that provides an estimate of the Hausdorff measure of a purely $k$-unrectifiable set intersected with suitable ``intrinsic tubes" \eqref{eq.assumpt.2}.
The consequence is a density upper bound estimate \eqref{eq:densityest},
that is the final step to conclude the proof of the implication.

Rather surprisingly, the classical approach can be made to work also in noncommutative homogeneous groups, see for instance \cite[Lemma 3.3.6]{Federer69} and \cite[Theorem 15.19]{Mat95}. However, some technical difficulties have to be overcome. We have introduced a definition of distance in the horizontal Grassmannian, that relies on horizontal projections and in a rather ``natural way"
leads to the ``distance estimates" of Theorem~\ref{thm:distest}.
It is of crucial importance the fact that the geometric constants of these estimates are independent of the horizontal subgroups.

The paper is structured as follows. 
Section~\ref{basic} recalls the basic definitions on homogeneous groups and introduces the horizontal Grassmannian. Section~\ref{technical} is devoted to definitions and preliminary facts that will be used to prove our results. We introduce tangent measures, the $(k,\G)$-approximate tangent group, differentiability results, the metric Jacobian, an area formula and a characterization of Radon measures with unique tangents, due to Mattila, \cite{Mat05}.
In Section~\ref{density} a density estimate on purely $k$-unrectifiable sets is established.
Section~\ref{main} contains the proof of our main result.

\section{Basic facts on homogeneous groups}\label{basic}

%%%%%%%%%%%%%%%%
\subsection{Graded nilpotent Lie groups and their metric structure}\label{sect:Graded}
A connected and simply connected {\em graded nilpotent Lie group} is a graded linear space $\G=H^1\oplus\cdots\oplus H^\iota$ equipped with a polynomial group operation such that its Lie algebra $\Lie(\G)$ is {\em graded}. This grading corresponds to the following conditions
\beq\label{eq:LieG}
\cg := \Lie(\G)=\cV_1\oplus\cdots\oplus\cV_\iota, \qquad [\cV_i,\cV_j]\subset \cV_{i+j}
\eeq
for all integers $i,j\ge0$ and $\cV_{j}=\{0\}$ for all $j>\iota$, with $\cV_\iota\neq\{0\}$. 
The integer $\iota \ge 1$ is the {\em step} of the group.
We denote by $\q$ the topological dimension of $\G$.

The graded structure of $\G$ agrees with a one parameter family of intrinsic dilations 
\[
\delta_r:\G\to\G
\]
defined as linear mappings such that
$\delta_r(p)=r^ip$ for each $p\in H^i$, $r>0$ and $i=1,\ldots,\iota$.
The graded nilpotent Lie group $\G$ equipped with intrinsic dilations
is called {\em homogeneous group}, \cite{FS82}.
With the stronger assumption that 
\beq
[\cV_1,\cV_j]=\cV_{j+1}
\eeq
for each $j=1,\ldots,\iota$ and $[\cV_1,\cV_\iota]=\set{0}$, we say that $\G$ is a 
{\em stratified group}. Further identifying $\G$ with the tangent space $T_0\G$ at the origin $0$, we have a canonical isomorphism between $H^j$ and $\cV_j$, that associates to each $v\in H^j$ the unique left invariant vector field $X\in\cV_j$ such that $X(0)=v$.

A Lie product defined on $\G$ induces a Lie algebra structure, where its group operation is given by the Baker-Campbell-Hausdorff (BCH) formula:
\beq\label{eq:BCH}
xy=\sum_{j=1}^\iota c_j(x,y)=x+y+\frac{[x,y]}2+\sum_{j=3}^\iota c_j(x,y)
\eeq
with $x,y\in\G$. Here $c_j$ denote homogeneous polynomials of degree $j$ with respect to the nonassociative Lie product on $\G$.

It is always possible to have these additional conditions, since the exponential mapping 
\[
\exp:\Lie(\G)\to\G
\]
of any simply connected nilpotent Lie group $\G$ is a bianalytic diffeomorphism.
In addition, the given Lie product and the Lie algebra associated to
the induced group operation are compatible, according to the 
following standard fact.

\bpr\label{pr:GAlgGroupOpe}
Let $\G$ be a nilpotent, connected and simply connected Lie group
and consider the new group operation given by \eqref{eq:BCH}.
Then the Lie algebra associated to this Lie group structure
is isomorphic to the Lie algebra of $\G$.
\epr
A {\em homogeneous distance} $d$ on a graded nilpotent Lie group $\G$ is a left invariant distance with $d(\delta_rx,\delta_ry)=r\,d(x,y)$ for all $x,y\in\G$ and $r>0$. 
We define \[ B(p,r)=\lb q\in\G: d(q,p)<r\rb\quad\mbox{and}\quad \B(p,r)=\lb q\in\G: d(q,p)\le r \rb \]
to be the open and closed balls, respectively. The corresponding homogeneous norm is denoted by 
\[
\|x\|=d(x,0)\quad\mbox{ for all $x\in\G$.}
\] 
The homogeneous distance $d$, along with its associated homogeneous norm $\|\cdot\|$, will be fixed for the sequel. When the graded nilpotent Lie group is equipped with the corresponding dilations, along with
a homogeneous norm, it is called {\em homogeneous group}.

Throughout the paper the symbol $\G$ will denote
a homogeneous group, equipped with a homogeneous distance $d$, if not otherwise stated.
The linear structure allows us to fix a scalar product
that makes $\G$ a finite dimensional Hilbert space. 
We also assume that all the layers $H^i$ are orthogonal to each other, namely the scalar product is {\em graded}. The symbol $|\cdot|$ is used to denote the Euclidean norm associated to this scalar product, that will be understood in the sequel.

\begin{Def}\label{d:hsubgroup}\rm
A linear subspace $S$ of $\G$ that satisfies $\delta_r(S)\subset S$ for every $r>0$ is a {\em homogeneous subspace} of $\G$.
If in addition $S$ is a Lie subgroup of $\G$ then we say that
$S$ is a {\em homogeneous subgroup} of $\G$.
\ed
Using dilations it is not difficult to check that $S\subset\G$ is a 
homogeneous subspace if and only if we have the direct decomposition
\[
S=S_1\oplus\cdots\oplus S_\iota,
\]
where each $S_j$ is a subspace of $H^j$.

\begin{Def}[Graded basis] \rm A {\em graded basis} $(e_1,\ldots,e_\q)$ of a homogeneous group $\G$ is a basis of vectors such that 
\beq
(e_{\m_{j-1}+1},e_{\m_{j-1}+2},\ldots,e_{\m_j})
\eeq
is a basis of $H^j$ for each $j=1,\ldots,\iota$, where we have set
\beq
\m_j=\sum_{i=1}^j\h_i \qandq \h_j=\dim H^j,
\eeq
with $\m_0 =0$. We will also denote $\m=\m_1$.
\ed

A graded basis provides the associated {\em graded coordinates} $x=(x_1,\ldots,x_\q)\in\R^\q$, that define the unique element $p=\sum_{j=1}^\q x_je_j\in\G$.

\br
It is easy to realize that one can always equip a homogeneous subgroup with graded coordinates.
\er

We conclude this section by recalling the definition of Hausdorff measure and
of spherical measure in metric spaces.

\bd
Let $X$ be a metric space with distance $d$ and let $A\subset X$, $k\in (0,\infty)$ and $\delta \in (0, \infty)$. The \textit{$k$-dimensional Hausdorff measure} $\cH_d^k$ is defined as
\[
\cH_d^k(A):=\sup_{\delta>0}\cH_{d,\delta}^k(A),
\]
where $\cH_{d,\delta}^k(A)=\inf \left\{ \sum_i 2^{-k}\diam_d(E_i)^k: A\subset \bigcup_i E_i, \diam_d(E_i)\le \delta\right\}$.
The spherical Hausdorff measure $\cS_d^k$ is obtained by requiring that the countable family $\set{E_i}$ is made only by closed metric balls. 

When the metric space is a homogeneous group equipped with a homogeneous distance we denote the $k$-dimensional Hausdorff measure by $\cH^k$. If the Hausdorff measure is considered in the Euclidean space equipped with the Euclidean norm $|\cdot|$, we denote it by $\cH^k_{|\cdot|}$.
\ed

%%%%%%%%%%%%%%%%
\subsection{The horizontal Grassmannian}\label{sect:intrinsicGrassmannian}

Let $S$ be a homogeneous subgroup of $\G$ as in Definition \ref{d:hsubgroup} and let
\[
S=S_1\oplus\cdots\oplus S_\iota
\]
be the decomposition induced by the layers of $\G$. We say that $S$ is a {\em horizontal subgroup} if $ S = S_1$ and a {\em vertical subgroup} if $S_j = H^j$ for all $2 \leq j \leq \iota$.

\br
It is easy to realize that any vertical subgroup is also
a normal subgroup.
\er

\bd
We say that two homogeneous subgroups $\V$ and $\W$ of $\G$ are {\em complementary subgroups} in $\G$ if $ \V \cap \W = \{0\}$ and $\G = \W \V$. If, in addition, $\W$ is normal we say that $\G$ is the {\em semidirect} product of $\V$ and $\W$ and write $\G = \W \rtimes \V$.
\ed

If $\G$ is the semidirect product of homogeneous subgroups $\V$ and $\W$, then we can define unique projections $\pi_\V : \G \to \V$ and $\pi_\W : \G \to \W$ in such a way that
\beq\label{eq:factorId}
\mathrm{id}_\G = \pi_\W \pi_\V.
\eeq
Furthermore, if $\W$ is normal in $\G$, then the following algebraic equalities hold:
\begin{equation} \label{eq.lemma.proj.alge} \begin{aligned}
& \pi_\W(p^{-1}) = \pi_\V(p)^{-1} \pi_\W(p)^{-1} \pi_\V(p), \quad \pi_\V(p^{-1}) = \pi_\V(p)^{-1}, \\
& \pi_\W(\delta_r p) = \delta_r \pi_\W(p), \quad \pi_\V(\delta_r p) = \delta_r \pi_\V(p), \\
& \pi_\W(p q) = \pi_\W(p) \pi_\V(p) \pi_\W(q) \pi_\V(p)^{-1}, \quad \pi_\V(p q) = \pi_\V(p) \pi_\V(q),
\end{aligned} \end{equation}
for $p,q\in\G$ and $r>0$. 

\bd\label{d:hlinMap}
Let $\P$ and $\G$ be two homogeneous groups. We say that a linear mapping $L:\P\to\G$ is an {\em h-homomorphism} if  
\[
L(xy)=L(x)L(y)\quad\text{and}\quad L(\delta^\P_rx)=\delta^\G_rL(x)
\]
for all $x,y\in\P$ and $r>0$,
where $\delta_r^\G$ and $\delta_r^\P$ represent
the dilations on $\G$ and $\P$, respectively.
\ed

\begin{Rem}\label{r:hlin}\rm
From \eqref{eq.lemma.proj.alge}, it follows that $\pi_\V$ 
is an h-homomorphism.
\end{Rem}

\begin{Def} \rm 
Let $\V$ be a horizontal homogeneous subgroup of $\G$. We denote by $\V^\perp$ the {\em orthogonal} homogeneous subgroup
\[
\V^\perp := V_1^\perp \oplus H^2 \oplus \cdots \oplus H^\iota,
\]
where $V_1^\perp$ is meant as the orthogonal of $V_1 = \V$ in $H^1$, with respect to fixed scalar product.
\ed
It is obvious to observe that $\V^\bot$ is a vertical subgroup.
%
% GRASSMANNIANS
%
\begin{Def}\label{d:abelian dimension}\rm
To every homogeneous group $\G$, we can associate a positive integer $\upsilon \leq \q$, which is the maximal linear dimension among all horizontal subgroups of $\G$.
\end{Def}

\bd \label{def:grass}
Let $\V$ be a $k$-homogeneous subgroup of $\G$. We say that $\V$ belongs to the 
family of subspaces $\cH(\G, k)$ if it is a $k$-dimensional and horizontal
subgroup. We say that $\cH(\G,k)$ is the {\em horizontal Grassmannian of $k$-dimensional subspaces}.
We write
\[
\cH(\G)=\bigcup_{k=1}^\upsilon\cH(\G, k)
\]
to indicate the {\em horizontal Grassmannian of $\G$}.
\ed

Notice that the horizontal Grassmannian $\cH(\G)$ is a compact subset of
the (standard) Grassmannian $\cG(\G)$, that is made by all linear subspaces of $\G$ and using our fixed Euclidean distance on $\G$.
Using the homogeneous distance $d$, we define the distance between two elements $\V_1,\V_2\in\cH(\G,k)$
as follows
\[
\rho(\V_1, \V_2)=\max_{\|x\|=1}d(\pi_{\V_1}(x),\pi_{\V_2}(x)),
\]
where $\pi_{\V_1},\pi_{\V_2}$ are the orthogonal projections onto $\V_1$ and $\V_2$, respectively. 
This distance defines the same topology of the standard distance on Grassmannians, but it turns out
to be more manageable in the setting of homogeneous groups.
In particular, we point out that $\cH(\G,k)$ is compact with respect to this topology.

In the sequel, we will frequently use the distance function
\[
d(p,S)=\inf\set{d(p,s):s\in S}=\inf\set{\|p^{-1}s\|:s\in S}
\]
with respect to our fixed homogeneous distance $d$.

\begin{The}\label{thm:distest}
For any $\V\in\cH(\G)$ the semidirect factorization 
\[
\G = \V^\perp \rtimes \V
\]
holds. In addition, there exists a universal constant $c_\G\in(0,1)$ such that 
for every $p\in\G$ we have
\begin{equation}\label{eq.piVcontin}
\| \pi_\V(p) \|\le c_\G^{-1} \|p\|
\end{equation}
and the following estimates hold
\begin{equation} \label{eq.lemma.proj3}
\begin{aligned}
& c_\G \, \|\pi_\V(p)\| \leq d(p, \V^\perp) \leq \|\pi_\V(p)\|,
\\ & c_\G\, \|\pi_\V(p)^{-1}  \pi_{\V^\perp}(p)  \pi_\V(p)\| \leq d(p, \V) \leq \|\pi_\V(p)^{-1} \pi_{\V^\perp}(p) \pi_\V(p)\|. \end{aligned}
\end{equation}
\et

\begin{proof}
We consider the orthogonal projection $P_{V_1} : H^1 \to V_1$  and set
\[
\pi_\V : \G \to \V, \quad \pi_\V(h_1 + \cdots +  h_\iota) := P_{V_1}(h_1),
\]
where $h_i\in H^i$. The we define 
\[
\pi_{\V^\bot}(p) = p\pi_\V(p)^{-1}=(p-\pi_\V(p))+
\sum_{j=2}^\iota c_j(p,-\pi_\V(p)).
\]
The first addend belongs to $V_1^\bot\subset H^1$, while \eqref{eq:BCH} joined with the grading assumption \eqref{eq:LieG}
shows that the second addend is in $H^2\oplus\cdots\oplus H^\iota$. We have established that the image of $\pi_{\V^\bot}$ 
is contained in $\V^\bot$. By definition of $\V^\bot$ we have
$\V\cap \V^\bot=\set{0}$. The definition of $\pi_{\V^\bot}$ yields \eqref{eq:factorId} with $\W=\V^\bot$,
therefore $\V^\bot\V=\G$. 
Due to the homogeneity properties of \eqref{eq.lemma.proj.alge}, for any $p\neq0$ we clearly have 
\begin{equation} \label{eq.1stEstim}
\frac{\| \pi_\V(p) \| + \| \pi_{\V^\bot}(p)\|}{\|p\|}\le \max_{\|z\|=1}\pa{\| \pi_\V(z) \| + \| \pi_{\V^\bot}(z)\|}.
\end{equation}
The crucial point is to obtain a universal bound on the right hand side, that is independent of $\V$.
We fix an arbitrary $\V_0\in\cH(\G,k)$ and observe that
\[
\max_{\|z\|=1}\| \pi_\V(z) \|\le \max_{1\le k\le\upsilon}\max_{\V\in\cH(\G,k)}d(\V,\V_0)+\max_{\|z\|=1}\|\pi_{\V_0}(z)\|=c_1(\G).
\]
In addition, it holds
\[
\max_{\|z\|=1}\| \pi_{\V^\bot}(z)\|\le 1+c_1(\G).
\]
Setting $c_\G=1/[1+2c_1(\G)]$, by \eqref{eq.1stEstim}, we have proved that
\begin{equation} \label{eq.10232}
c_\G\pa{\| \pi_\V(p) \| + \| \pi_{\V^\bot}(p)\|}\le  \|p\|.
\end{equation}
In particular, \eqref{eq.piVcontin} holds.
To prove \eqref{eq.lemma.proj3}, we start by noticing that
\[
d(p, \V^\perp) = \inf \{ \|w^{-1} p \|  :  w \in \V^\perp \} \leq \|\pi_\V(p)\|,
\]
by plugging in $w = \pi_{\V^\perp}(p)$. In a similar fashion, we apply \eqref{eq.10232} to deduce that
\[ 
\| w^{-1} p \| = \| w^{-1} \pi_{\V^\perp}(p) \pi_\V(p) \| \geq c_\G( \|w^{-1} \pi_{\V^\perp}(p)\| + \|\pi_\V(p)\|) \geq c_\G \| \pi_\V(p) \|. 
\]
The distance $d(p, \V)$ can be bounded from below as before by taking $v = \pi_\V(p)$ in
\[
d(p, \V) = \inf\{ \| v^{-1} p \|  :  v \in \V \} \leq \| \pi_\V(p)^{-1} \pi_{\V^\perp}(p) \pi_\V(p) \|.
\]
For the opposite bound, we apply once again \eqref{eq.10232} and obtain
\[ \begin{aligned}
d(p, \V) & = \inf\{ \| p^{-1} v \|  : v \in \V \} =
\\ &= \inf\{ \| \pi_\V(p)^{-1} \pi_{\V^\perp}(p)^{-1} \pi_\V(p) \pi_\V(p)^{-1} v \|  : v \in \V \} \geq
\\ & \geq c_\G \left[ \inf\{ \| \pi_\V(p)^{-1} \pi_{\V^\perp}(p)^{-1} \pi_\V(p) \| + \| \pi_\V(p)^{-1} v \|  : v \in \V \} \right] =
\\& = c_\G \| \pi_\V(p)^{-1} \pi_{\V^\perp}(p)^{-1} \pi_\V(p) \|.
\end{aligned} \]
We observe that 
\[
\|\pi_\V(p)^{-1} \pi_{\V^\perp}(p)^{-1} \pi_\V(p) \|=\|(\pi_\V(p)^{-1} \pi_{\V^\perp}(p)^{-1} \pi_\V(p))^{-1} \|,
\]
therefore the previous estimates complete the proof.
\end{proof}

The universal constant $c_\G$ of the previous theorem will appear several times in the sequel. An immediate consequence of \eqref{eq.lemma.proj3} is the estimate
\begin{equation} \label{eq.2.18}
d(p, \pi_\V(p)) \leq \frac{1}{c_\G} d(p, \V)
\end{equation}
for all $p\in\G$ and $\V\in\cH(\G)$.

\section{Technical preliminaries}\label{technical}

This section is devoted to the essential tools that will be used to prove our main results.

\subsection{Tangent measures} For $p \in \G$ and $r > 0$, similarly to the Euclidean setting, we can define the {\em magnification map} using intrinsic dilations, by setting
\[
T_{p, r}(q) := \delta_{\frac{1}{r}}(p^{-1} q)
\]
for all $p,q\in\G$ and $r>0$.
Notice that $T_{p, r}$ maps the ball of center $p$ and radius $r$ to the unit ball centered at the identity element of $\G$.

\bd
Let $\mu$ be a Radon measure on $\G$ and let $p \in \G$. A Radon measure $\nu$ on $\G$ with $\nu(\G)>0$ is said to be a {\em tangent measure} of $\mu$ at $p$, if there are sequences of positive numbers $c_i$ and $r_i \to 0$ such that
the following weak convergence of measures holds
\[
c_i (T_{p, r_i})_\# \mu \rightharpoonup \nu
\]
as $i\to\infty$. The set of all tangent measures of $\mu$ at $p$ is denoted by $\Tan(\mu, p)$.
\ed
Notice that if $\nu \in \Tan(\mu, p)$, clearly $\lambda \nu \in \Tan(\mu, p)$ for all $\lambda > 0$, hence uniqueness is understood up to a positive factor.
The next theorem is a consequence of combining \cite[Theorem~3.2]{Mat05} and \cite[Lemma~2.5]{Mat05}.
It holds for general locally compact metric Lie groups equipped with dilations.

\bt \label{thm.tangentmeasures}
Let $\mu$ be a Radon measure on $\G$. The following conditions are equivalent: \mbox{}
\begin{enumerate}[label=\textbf{(\arabic*)}, itemsep=.6em]
\item At $\mu$-a.e. $p \in \G$ there exists a unique tangent measure $\nu_p \in \Tan(\mu, p)$.
\item For $\mu$-a.e. $p \in \G$ there exists a closed homogeneous subgroup $\V_p$ of $\G$ for which
\[
\Tan(\mu,p) = \{ \lambda \nu_p  :  0 < \lambda < \infty \},
\]
where $\nu_p$ is a Haar measure on $\V_p$.
\end{enumerate}
If one of these two conditions hold, then for $\mu$-a.e. $p\in\G$ there exists $c_p>0$ such that 
\[
\frac{(T_{p,r})_\#\mu}{\mu(B(p,r))}\to c_p \nu_p\quad\text{as}\quad r\to0^+.
\]
\et

\subsection{Density and approximate tangent group} Let $\mu$ be a Radon measure in $\G$. We define the {\em upper} and {\em lower $k$-densities} of $\mu$ at $p \in \G$ as
\[ \Theta^{\ast k}(\mu, p) := \limsup_{r \to 0^+} \frac{\mu(B(p, r))}{r^k} \quad \text{and} \quad \Theta_\ast^k(\mu, p) := \liminf_{r \to 0^+} \frac{\mu(B(p, r))}{r^k}.
\]
If $E \subset \G$ is $\cH^k$ measurable and $\mu = \cH^k \res E$ is locally finite, then we define the corresponding upper and lower $k$-densities 
\[
\Theta^{\ast k}(E, p)= \limsup_{r \to 0^+} \frac{\cH^k(E\cap B(p, r))}{r^k} \quad\text{ and } \quad \Theta_\ast^k(E, p)= \liminf_{r \to 0^+} \frac{\cH^k(E\cap B(p, r))}{r^k} ,
\]
respectively.
The next result recalls the standard density estimates for Hausdorff measures, see for instance \cite[2.10.19]{Federer69}. 

\bl \label{lemma:densities}
Let $E \subset \G$ be $\cH^k$-measurable with $\cH^k(E)<+\infty$. Then 
\begin{enumerate}[itemsep =0.4em, label=(\roman*)]
\item For $\cH^k$-a.e. $p \in E$,
\[ 2^{-k} \leq \Theta^{\ast k}(E, p) \leq 1. \]
\item For $\cH^k$-a.e. $p \in\G\setminus E$, $\Theta^{\ast k}(E, p) = 0$.
\end{enumerate}
\el

It is well known in the Euclidean setting that rectifiability can be characterized by cones, see \cite[Theorem~15.19]{Mat95}.
The next definition provides the proper notion of cone in homogeneous groups, that naturally works with the notion of rectifiability.

\begin{Def}\rm
Let $s \in (0,\,1)$, $p \in \G$ and let $\T$ be a homogeneous subgroup of $\G$. The {\em intrinsic cone} of vertex $p$, axis $\T$ and opening $s$ is defined as
\[
X(p,\,\T,\,s)= \left\{ q \in \G  : d(p^{-1} q,\,\T) < s d(p,\,q) \right\}.
\]
We will also use the notation $X(p,r,\T,s)=X(p,\T,s)\cap B(p,r)$.
\ed

The notion of intrinsic cone leads us to the notion of ``intrinsic approximate tangent group'' for an $\cH^k$-measurable set.
\bd
Let $E \subset \G$ be an $\cH^k$-measurable set. A homogeneous $k$-subgroup $\T_p$ of dimension $k$ and Hausdorff dimension $k$ is a {\em $(k,\G)$-approximate tangent group} to $E$ at $p$ if the following properties hold: \mbox{}
\begin{enumerate}
\item $\Theta^{\ast k}(E,p) > 0$,
\item for all $s \in (0,1)$ we have
\[ \lim_{r \to 0^+} \frac{ \cH^{k} \res E ( B(p,r) \setminus X(p,\T_p,s))}{r^{k}} = 0.  \]
\end{enumerate}
We denote by $\mathrm{ap}\Tan_\G^k(E,p)$ the set of all $(k,\G)$-approximate tangent groups to $E$ at $p$ and we use $\T_p$ when it is unique.
\ed

The next result establishes several properties of the $(k,\G)$-approximate tangent group, along with its uniqueness up to negligible sets.

\bpr \label{prop:uniqaptan}
Let $E \subset \G$ be an $\cH^k$-measurable set, $1 \leq k \leq \upsilon$, and let
\[
A = \{ p \in E:\mathrm{ap}\Tan_\G^k(E,p) \neq \emptyset \}.
\]
Then the following properties hold: \mbox{}
\begin{enumerate}[label=(\roman*), itemsep=0.4em]
    \item The set $A$ is $\cH^k$-measurable.
    \item For $\cH^k$-a.e. $p \in A$ there exists a unique $(k,\G)$-approximate tangent group $\T_p$.
    \item The map $p \mapsto \T_p$ is measurable.
\end{enumerate}
\epr
The proof of this theorem does not involve the type of algebraic structure of the group, rather it is measure theoretic. Indeed, it can be carried out following the same argument of \cite[Proposition 3.9]{MSSC10} for the Heisenberg group.

\subsection{Rectifiability, differentiability and metric Jacobian} Next we recall the standard notion of rectifiability in metric spaces, see for instance  \cite[3.2.14]{Federer69}.

\begin{Def}[Rectifiable set] \label{def.rectifiable} \rm A measurable set $E \subset \G$ is $k$-rectifiable if there exist sets $A_i\subset\R^k$ and Lipschitz maps
$f_i : A_i \subset \R^k \to \G $, with $i\in\N$, such that
\[
\cH^k\Big(E\sm\bcup_{i\in\N}f_i(A_i)\Big)=0.
\]
\ed
\br\label{rmk.closed} %%
By completeness of $\G$, we may assume that all $A_i\subset\R^k$ in the definition of $k$-rectifiability are closed sets. 
\er

%Under these conditions, up to recursively considering the differences $\widetilde{A}_i=A_i\sm f_i^{-1}\left(\bigcup_{j=0}^{i-1}f_j(\widetilde{A}_j)\right)$, we may assume $A_i$ to be Borel sets and all the sets $f_i(A_i)$ disjoint.

For Banach homogeneous group targets, the a.e.\ differentiability of Lipschitz mappings still holds if the first layer of the target satisfies the so-called Radon--Nikodym property, see \cite[Theorem~1.1]{MagnaniRajala2014}.
In particular, it holds for any Lipschitz mapping from a subset of a stratified group to a finite dimensional homogeneous group. 

\begin{Def}[Differentiability]\rm 
If $A\subset\R^k$, $f:A\to\G$, $x\in A$ is a density point and $L:\R^k\to\G$ is an h-homomorphism, we say that $f$ is differentiable at $x$ if 
\[ d(f(x)^{-1}f(xz),L(z))=o(d(z,0)) \]
as $z\to 0$. The differential of $f$ at $x$ is denoted by $Df(x)$.
\ed
The differential for group-valued mappings also inherits a Lie group homomorphism property, namely it is an
h-homomorphism, according to Definition~\ref{d:hlinMap}.

Differentiability is strictly related to area through the area formula. 
Kirchheim observed that to compute the area of a set parametrized by a Lipschitz mapping 
on a subset of Euclidean space a weaker notion of differentiability suffices, \cite{Kir94}. 
This is the so-called the ``metric differentiability''.

In analogy with the Euclidean framework, the metric differential on homogeneous groups is represented by a {\em homogeneous seminorm} $s:\G\to[0,+\infty)$, that is a continuous mapping satisfying the properties 
\[
s(\delta_rx)=rs(x)\quad \text{ and }\quad s(xy)\le s(x)+s(y)
\]
for all $x,y\in\G$ and $r>0$.

\begin{Def}[Metric differential]\rm
Let $A \subset \R^k$, $f : A \to \G$ and fix a density point $x \in A$. We say that $f$ is \textit{metrically differentiable} at $x$ if there exists a homogeneous seminorm $s$ such that
\[ d( f(x), f(xz)) - s(z) = o(\|z\|) \quad \text{as $z \in x^{-1}A$ and $\|z\| \to 0^+$}. \]
The homogeneous seminorm $s$ is unique and it is denoted by $mdf(x)$, which we call the {\em metric differential of $f$ at $x$}.
\ed
\br \label{rmk.diff-metricdiff}
Differentiability implies metric differentiability. It is easy to notice that when $f:A\to\G$ is differentiable at $x$, then it is also metrically differentiable at $x$ and
\[
mdf(x)(h)= d(D f(x)(h),0) \quad \text{for $h \in \R^k$}.
\]
\er
A general notion of metric Jacobian can be introduced for metric space targets, see \cite{MagnaniRajala2014} for more information. 
\begin{Def}[Metric Jacobian]\rm
Let $f : A \subset \R^k \to \G$ be metrically differentiable at $x$ and let $s$ denote its metric differential $mdf(x)$.
The {\em metric Jacobian} is defined as follows
\[  Jf(x) = \begin{cases} \displaystyle\frac{ \cH_s^k(B_E(0,1)) }{ \cH^k(B_E(0,1)) } & \text{if $s$ is a homogeneous norm,}
\\[.9em] 0 & \text{otherwise}. \end{cases} \]
The unit ball $B_E(0,1)\subset\R^k$ is defined with respect to the Euclidean metric of $\R^k$.
\ed
With this notion, a general metric area formula holds for 
a.e.\ metric differentiable mappings defined on a stratified group.
More specifically, the following theorem is a consequence of combining \cite[Theorem~1.1]{MagnaniRajala2014} and \cite[Theorem~1.2]{MagnaniRajala2014}, when the Lipschitz mapping is defined
on a subset of the Euclidean space. 

\bt[\cite{MagnaniRajala2014}]\label{t:area}
Let $A$ be a measurable set of $\R^k$ and let $f:A\to\G$ be Lipschitz, where $\G$ is a homogeneous group. Then
\[
\int_A Jf(x) dx =\int_\G \sharp(f^{-1}(y))d\cH^k(y).
\]
\et
In our setting, a ``linearization type theorem'' is available. Previous versions of this result in more specific assumptions are \cite[Lemma~3.22]{Federer69}, \cite[Lemma~4]{Kir94} and \cite[Proposition 4.1]{Mag}.
Specializing \cite[Lemma~4.2]{MagnaniRajala2014} to our setting
and taking into account formula (51) of \cite{MagnaniRajala2014}, we obtain
the following result.

\bl\cite[Lemma~2.4]{MagnaniRajala2014}\label{l:linearization}
Let $A\subset\R^k$ be a closed set, let $f:A\to\G$ be Lipschitz
and let $\cD\subset A$ be the subset of differentiability points,
where the metric differential of $f$ is a homogeneous norm.
Then the following statements hold. 
\begin{enumerate}
\item
There exists a family of Borel sets $\{A_i\}_{i\in\N}$ such that 
$\cD=\bcup_{i\in \N} A_i$ and $f|_{A_i}$ is bi-Lipschitz onto its image for all $i\in\N$.
\item
For a.e.\ $x\in \cD$, we have 
\[
Jf(x)=\frac{ \cH^k\left(Df(x)(\B(0,1)\right) }{ \cH_{|\cdot|}^k(\B_{E}(0,1)) }=\limsup_{r\to0^+}\frac{\cH^k(f(E\cap \B(x,r)))}{\cH^k_{|\cdot|}(\B_{E}(x,r))}. 
\]
\end{enumerate}
\el
We have introduced the notation $\B_E(0,1)$ to denote the Euclidean unit ball of $\R^k$.
The following "bi-Lipschitz linearization" essentially extends the previous lemma 
to the case where the domain of the Lipschitz mapping is measurable.
Slightly more information is given, since we also consider the metric Jacobian. 
\bl \label{lemma.3.1} Let $E \subset \G$ be a $k$-rectifiable set. Then there exist
measurable set $M_i\subset\R^k$ and Lipschitz maps $f_i : M_i \to \G$ such that for every $i\in\N$ the following
conditions hold.
\begin{enumerate}
\item
$f_i$ is bi-Lipschitz onto its image
\item
$f_i$ is everywhere metrically differentiable 
\item
$mdf_i(x)$ is a homogeneous norm for all $x\in M_i$,
\[
Jf_i(x)=\frac{ \cH^k\left(Df_i(x)(\B(0,1))\right) }{\cH_{|\cdot|}^k(\B_E(0,1)) }
\]
for a.e.\ $x\in M_i$ and $x\to Jf_i(x)$ is measurable.
\end{enumerate}
\item
Finally, there exists an $\cH^k$-negligible set $E_0\subset\G$ such that
\beq\label{eq:inclusionE_0}
E \subseteq \bigcup_{i \in \N} f_i(M_i) \cup E_0
\eeq
and measurable sets $M_i$ can be chosen to make the family $\set{f_i(M_i):i\in\N}$ disjoint.
\el
\begin{proof}
By Remark~\ref{rmk.closed}, \cite[Theorem~1.1]{MagnaniRajala2014}, 
Lemma~\ref{l:linearization} and Theorem~\ref{t:area}, we end up with an $\cH^k$-negligible set $E_0$ and a countable family of Borel sets $B_i\subset\R^k$ and Lipschitz mappings $f_i:B_i\to\G$ that satisfy condition (1), (2) and (3)
of our claim, along with \eqref{eq:inclusionE_0}. Taking into account that $f_i$ are bi-Lipschitz, then their images and preimages of Borel sets are $\cH^k$ measurable and $\cL^k$ measurable, respectively. Thus, the measurable sets 
\[
M_i=B_i\sm f_i^{-1}\pa{\bcup_{l=1}^{i-1}f_l(M_l)}
\]
recursively defined makes the images $f(M_i)$ disjoint. 
\end{proof}

\section{Density estimate on purely unrectifiable sets}\label{density}

The main result of this section is Theorem~\ref{ppwpod}, that can be seen as
an improving of a decay estimate on suitable intersections of cones with 
purely $k$-unrectifible sets.

\bd
We say that a set $E\subset\G$ is {\em purely $k$-unrectifiable} if for every Lipschitz mapping $f:A\to\G$ with $A\subset\R^k$ we have
\[
\cH^k(E\cap f(A))=0.
\]
\ed
Notice that in \cite[3.2.14]{Federer69} the set of the previous definition would have been called purely $(\cH^k,k)$-unrectifiable.
\bl \label{lemma.5.1}
If $E \subset \G$, $\T \in \cH(\G,k)$, $r > 0$, $0<s<1$, and
\beq\label{eq:EcapXpTperp}
E \cap X(p,\T^\perp,s)\cap B(p,r) = \emptyset
\eeq
whenever $p \in E$, then every subset of $E$ with diameter less than $r$ is contained in the image of some Lipschitz map $f : \R^k \to \G$ with $\mathrm{Lip}(f) \leq s^{-1}$.
\el

\begin{proof}
Let $F \subset E$ be a subset with $\mathrm{diam}(F)<r$ and let $\pi_\T$ be the orthogonal projection onto $\T$. If $p,q \in F$, then $\|p^{-1} q\| < r$ and by our assumption $q \notin X(p,r,\T^\bot),s)$.
It follows that
\[
s^{-1} d( p^{-1} q, \T^\perp) \geq \|p^{-1} q\|.
\]
The estimates \eqref{eq.lemma.proj3} on the distance with respect to vertical subgroups leads to the inequality
\[
s^{-1} \| \pi_{\T}(p^{-1} q) \| \geq \|p^{-1} q\|.
\]
Finally, the homomorphism property of $\pi_\T$ immediately gives our claim.
\end{proof}

\begin{Rem}\label{r:krect}\rm
In the assumptions of the previous lemma, it is easy to realize that
replacing the condition \eqref{eq:EcapXpTperp} by 
\[
E \cap X(p,\T^\perp,s) = \emptyset,
\]
we then obtain a unique Lipschitz mapping that parametrizes $E$.
\end{Rem}

\bt\label{ppwpod}
Let $E\subset\G$ be a purely $k$-unrectifiable set, consider $\V \in \cH(\G,k)$, $0 < s < c_\G^3$, $\lambda > 0$, and $\delta>0$, where $c_\G$ is 
as in \eqref{eq.lemma.proj3}. If for all $p \in E$ and any $0 < r \leq \delta$ we have
\begin{equation}\label{eq.assumpt.1}
\cH^k \res E \left( X(p,r,\V^\perp,s)\right) \leq \lambda r^k s^k,
\end{equation}
then for all $w\in\G$ and $0 < t \le \delta/6$ we get
\begin{equation}\label{eq.assumpt.2}
\cH^k \res E \left( B \left(w,t \right) \cap \pi_{\V}^{-1}\pa{B \left (\pi_{\V}(w),\rho} \right) \right) \leq 2 \,\lambda\, (21)^k\, t^k.
\end{equation}
In particular, for all $w\in\G$ we have the density upper bound
\begin{equation}\label{eq:densityest}
\Theta^{\ast k}(\cH^k\res E,w) \leq 2\, (21)^k\, \lambda.
\end{equation}
\et

\begin{proof}
We fix $w \in \G$ and $0 < \rho \leq s\delta/(24\,c_\G)$.
We set the parameter $\epsilon=s/4c_\G$ and define
\[
A = E \cap B\left(w, \frac{\rho}{\epsilon} \right) \cap \pi_{\V}^{-1}\pa{ B\left(\pi_{\V}(w),\rho\right)}.
\]
We consider the set
$
C =\{y \in A : A \cap X(y, \V^\perp,\epsilon) \neq \emptyset \}.
$
%more details in details-RG.lyx
Taking into account Lemma \ref{lemma.5.1} and Remark~\ref{r:krect}, the difference $A\setminus C$ is $k$-rectifiable, therefore $\cH^k(A\setminus C) = 0$
and we can prove \eqref{eq.assumpt.2} replacing $A$ by $C$.
For every $x \in C$, the following function is well defined
\[
h(x)= \sup \left\{ \|x^{-1} y\|  : y \in A \cap X(x,\V^\perp,\epsilon) \right\}
\]
and it satisfies the easy estimates $0 < h(x) \leq \diam (A) \leq 2 \rho/\epsilon$. We have the inclusion
\begin{equation}\label{eq.inc.1}
\pi_\V \left( C \right) \subset \bigcup_{x \in C} B\left(\pi_\V(x), \frac{\epsilon h(x)}{5}\right) \cap \V \subset B\left( \pi_\V(w), \frac{7}{5} \rho\right)\cap \V,
\end{equation}
since for all points $x$ of $C$ we have $d(\pi_\V(x),\pi_\V(w))<\rho$.
We apply Vitali's covering lemma to find an at most countable subset $D \subset C$ such that
\[
\left\{ B\left(\pi_\V(x), \frac{\epsilon h(x)}{5}\right)  : x \in D \right\}
\]
is disjoint and the corresponding family of balls with radii $\epsilon h(x)$ covers $\pi_\V(C)$. We have
\begin{equation}\label{eq:inclusionD}
\pi_\V(C) \subseteq \bigcup_{x \in D} B(\pi_\V(x), \epsilon h(x)) \cap \V.
\end{equation}
In particular, $\pi_\V$ restricted to $D$ is injective. 
Since $\V$ is a homogeneous subgroup of Hausdorff dimension $k$, there exists a constant
$\gamma_\V>0$ such that 
\[
\cH^k=\gamma_\V\cL^k,
\]
where $\cL^k$ is the $k$-dimensional Lebesgue measure on $\V$ with respect to an
orthonormal system of coordinates on $\V$. It follows that
\[
\cH^k(\V\cap B(u,t))=\gamma_\V\cL^k(\V\cap B(0,1))\,t^k
\]
for every $u\in \V$ and $t>0$. 
As a result, the inclusion \eqref{eq.inc.1} immediately implies that
\begin{equation}\label{eq.5.111}
\sum_{x \in D} \epsilon^k h(x)^k \leq 7^k \rho^k.
\end{equation}
To obtain our claim we have to estimate the measure
\[
\cH^k \res C \left( \pi_{\V}^{-1}(B(\pi_{\V}(x),\epsilon h(x))) \right) 
\]
for all $x \in D$ by $\lambda \epsilon^k\, h(x)^k$, up to a universal geometric factor,
and then use \eqref{eq.5.111}.
Let us fix $x \in D$. By definition of $h(x)$, we may find 
$y \in A \cap X(x,\V^\perp,\epsilon)$ such that
\beq\label{eq:43rd}
4 \|x^{-1} y\| > 3 h(x).
\eeq
Using \eqref{eq.lemma.proj3}, we immediately notice that
\[
c_\G \|\pi_\V(x^{-1} y)\| < \epsilon \|x^{-1} y\| \leq \epsilon h(x),
\]
therefore $\|\pi_\V(x^{-1} y)\| < \epsilon h(x)/c_\G$.
We claim that $C \cap \pi_\V^{-1}(B(\pi_\V(x), \epsilon h(x)))$ is contained in
the following union of cones 
\[
X\Big(x, 3h(x),\V^\perp,\frac{s}{c_\G^2}\Big) \cup X\Big(y,3h(x),\V^\perp,\frac{s}{c_\G^2}\Big).
\]
Let $z \in C$ be such that $\|\pi_\V(x^{-1} z)\| < \epsilon h(x)$ and notice that either 
$z\in A \cap X(x,\V^\perp,\epsilon)$, therefore $\|x^{-1} z\| \leq h(x)$ or
otherwise $z \notin X(x,\V^\perp,\epsilon)$, hence
\[
\|x^{-1} z\| \leq \frac{1}{\epsilon} \|\pi_\V(x^{-1} z)\| < h(x).
\]
Therefore, in both cases $\| x^{-1} z \| \leq h(x)<3h(x)$ and consequently
\[
\| y^{-1} z \| \leq \| x^{-1} z \| + \| x^{-1} y\| \leq 2 h(x)<3h(x).
\]
It follows that $z\in B(x,3h(x))\cap B(y,3h(x))$ and 
\[ \begin{aligned}
\| \pi_\V(x^{-1} z) \| + \| \pi_\V(y^{-1} z) \| & \leq 2 \| \pi_\V(x^{-1} z) \| + \| \pi_\V(x^{-1} y) \|
\\ & < 2\epsilon h(x)+\frac{\epsilon}{c_\G}h(x)\le\frac{3\epsilon}{c_\G} h(x).     
\end{aligned}
\]
Taking into account \eqref{eq:43rd}, we get
\[
\begin{aligned}
\| \pi_\V(x^{-1} z) \| + \| \pi_\V(y^{-1} z) \| & < \frac{4\epsilon}{c_\G} \|x^{-1}y\|=\frac{s}{c_\G^2}\|x^{-1}y\|   \\
&\leq \frac{s}{c_\G^2} \| x^{-1} z \| + \frac{s}{c_\G^2} \| x^{-1} y \|.
\end{aligned} \]
It follows that  
\[
C \cap \pi_\V^{-1}(B(\pi_\V(x), \epsilon h(x))) \subset E\cap 
\left(X\Big(x, 3h(x),\V^\perp,\frac{s}{c_\G^2}\Big) \cup X\Big(y,3h(x),\V^\perp,\frac{s}{c_\G^2}\Big)\right).
\]
Since $3h(x)\le 6\rho/\ep=24c_\G\rho/s\le \delta$, 
we may apply \eqref{eq.assumpt.1}, getting
\[
\cH^k \left(C \cap \pi_\V^{-1} (B(\pi_\V(x), \epsilon h(x)))\right) \leq 2 \lambda  \left(3h(x) \frac{s}{c_\G^2}\right)^k=
2\lambda\left(\frac{12\,h(x)\ep}{c_\G}\right)^k.
\]
As a consequence of \eqref{eq:inclusionD}, we observe that
\[
C\subseteq \bigcup_{x \in D} \pi_\V^{-1}\left(B(\pi_\V(x), \epsilon h(x))\right).
\]
Summing the previous estimate over all $x \in D$, using \eqref{eq.5.111} and taking into account that $\cH^k(A\sm C)=0$, we obtain
\[
\cH^k \res E \left( B \Big(w,\frac{4c_\G\rho}{s} \Big) \cap \pi_{\V}^{-1}\pa{B \left (\pi_{\V}(w),\rho} \right) \right) \leq 2 \left(\frac{84}{c_\G}\right)^k \lambda  \rho^k.
\]
Setting $t=4c_\G\rho/s$ the proof is concluded.
\end{proof}

\section{Proof of the main result}\label{main}

This section is devoted to the proof of Theorem~\ref{thm.principale}.

We start with the proof of ({\romannumeral 1}) $\implies$ ({\romannumeral 2}). In view of Lemma \ref{lemma.3.1},
we can find a countable family of bi-Lipschitz mappings $f_i:E_i\to S_i$, $E_i\subset\R^k$ is measurable, $S_i\subset\G$ is $\cH^k$ measurable, $i\in\N$, such that
\begin{enumerate}[label=(\alph*), itemsep=.5em]
\item the sets $S_i$ are disjoint,
\item $E$ is equal to the union of an $\cH^k$-negligible set $E_0$ with $\bigcup_{i\in\N} S_i$,
\item $mdf_i(x)$ is a homogeneous norm at all $x \in E_i$.
\end{enumerate}
The proof of \eqref{eq.blowup} for $\cH^k$-a.e.\ $p\in E$ is actually equivalent to the validity for all $i\in\N$ of the limit
\begin{equation} \label{eq.blowup2}
\frac{1}{r^k}(T_{p,r})_\#(\cH^k\res S_i) \rightharpoonup \cH^k \res D f_i(f_i^{-1}(p))(\R^k)
\end{equation}
at $\cH^k$-a.e.\ $p\in S_i$, being $Df_i(f_i^{-1}(p))(\R^k)\in\cH(\G,k)$. Indeed, using Lemma~\ref{lemma:densities}, for $\cH^k$-a.e.\ $p \in S_i$ we have
\begin{equation} \label{eq.blowup3}
\Theta^{\ast k}(E \setminus S_i,p) = 0.
\end{equation}
Then the two statements are equivalent, due to the equality 
\[ 
\frac{1}{r^k} \int_E \varphi \circ T_{p,r} \mathrm{d} \cH^k = \frac{1}{r^k} \int_{E \setminus S_i} \varphi \circ T_{p,r} \mathrm{d} \cH^k + \frac{1}{r^k} \int_{S_i} \varphi \circ T_{p,r} \mathrm{d} \cH^k
\]
for all $\varphi \in C_c(\G)$, where $p$ satisfies \eqref{eq.blowup3}.

For the sake of notation, in the sequel we denote $f_i$ by $f$, $M_i$ by $M$ and
$S_i$ by $S$. We may also fix $p=f(x)$ such that \eqref{eq.blowup3} holds, we have the estimates
\[ 
\begin{cases} \| Df(x)(y) \| \geq c_x |y| & \text{for all $y \in \R^k$}, \\
d(f(x), f(y)) \geq c |x-y| & \text{for all $y \in M$},
\end{cases} 
\]
with $c,c_x>0$. Here we have chosen $x\in M$ such that it is a Lebesgue point of $Jf$ and a differentiability poiny of $f$. 
Let $\varphi \in C_c(\G)$ and let $R > 0$ be such that $\supp\, \ph\subset B(0,R)$. Then
defining $C_x=\min\set{c,c_x}>0$, we get
\[ \varphi\left( \delta_{\frac{1}{r}}(f(x)^{-1}f(y))\right) = \varphi\left( D f(x) \left( \frac{y-x}{r} \right) \right) = 0 \]
whenever $|x - y| \geq \frac{Rr}{C_x}$. Since $x$ is a Lebesgue point for $Jf$
and $f$ is differentiable at $x$, we infer that
\beq\label{eq:o(r^k)}
\int_M \varphi\left( \delta_{\frac{1}{r}} (f(x)^{-1}f(y))\right) Jf(y) \mathrm{d}y - Jf(x) \int_M \varphi \left( D f(x)\left( \frac{y-x}{r}\right) \right)\mathrm{d}y= o(r^k)
\eeq
as $r \to 0^+$. From the area formula of Theorem~\ref{t:area}, it follows that
\beq \label{eq.main.1} \begin{aligned} \int_\G \varphi \, 
\mathrm{d}\cH^k \res D f(x)(\R^k) & = \int_{D f(x)(\R^k)} \varphi \, \mathrm{d}\cH^k 
\\ & = r^{-k} Jf(x) \int_{\R^k} \varphi \left( D f(x)\left( \frac{y-x}{r}\right) \right) \mathrm{d}y. \end{aligned} \eeq
The same area formula along with basic properties of push-forward measures yield
\beq \label{eq.main.2} \begin{aligned} 
\int_M \varphi\left( \delta_{\frac{1}{r}} (f(x)^{-1} f(y))\right) Jf(y) \mathrm{d}y & = \int_\G\varphi\left( \delta_{\frac{1}{r}} (p^{-1}z) \right) \mathrm{d}\cH^k\res S(z) 
\\ & = \int_{\G} \varphi \, \mathrm{d}(T_{p,r})_\symbol{35}( \cH^k \res S). \end{aligned}\eeq
Considering \eqref{eq.main.1}, \eqref{eq:o(r^k)} and \eqref{eq.main.2}, we find that
\[
\begin{aligned} \int_\G \varphi \, \mathrm{d}\cH^k \res Df(x)(\R^k)  & = \lim_{r \to 0^+} r^{-k} Jf(x) \int_{\R^k} \varphi \left( D f(x)\left( \frac{y-x}{r}\right) \right) \, \mathrm{d}y
\\ & = \lim_{r \to 0^+}r^{-k} \int_B \varphi
\left(\delta_{\frac1r}(f(x)^{-1} f(y)) \right) Jf(y) \, \mathrm{d}y  
\\& = \lim_{r \to 0^+} r^{-k} \int_{\G} \varphi \, \mathrm{d} \left( (T_{p,r})_\symbol{35}( \cH^k \res S) \right),\end{aligned} \]
hence concluding the proof of our claim.

Proof of ({\romannumeral 2}) $\implies$ ({\romannumeral 3}).
For $\cH^k$-a.e. $p\in E$, the assumption ({\romannumeral 2}) tells us that
\[
\frac{1}{r^k}(T_{p,r})_\# (\cH^k \res E) \rightharpoonup \cH^k \res \T_p
\]
as $r \to 0^+$. 
If $\nu_p \in \Tan(\cH^k \res E,p)$, then there are sequences $c_i\in\R$ and $r_i>0$ such that $r_i$ is infinitesimal and
\[
c_i (T_{p,r_i})_\# (\cH^k \res E) \rightharpoonup \nu_p.
\]
Therefore, being $\nu_p$ nonvanishing, the sequence $r_i^k c_i$ converges
to some $\lambda_p\in\R\sm\set{0}$ and the following holds:
\[
c_i (T_{p,r_i})_\# (\cH^k \res E) \rightharpoonup \lambda_p \, \cH^k \res \T_p. 
\]

Proof of ({\romannumeral 3}) $\implies$ ({\romannumeral 4}). By Theorem \ref{thm.tangentmeasures}, 
for a.e.\ $p\in E$ there exists $c_p>0$ such that
\[
\frac{1}{\cH^k \res E(B(p,r))}(T_{p,r})_\# (\cH^k \res E) \rightharpoonup c_p \cH^k \res \T_p
\]
as $r \to 0^+$. By homogeneity and the fact that $\cH^k\res \T_p$ is locally finite, we get
\[ \cH^k \big( \T_p \cap \partial( B(e,1) \setminus X(e,\T_p,s) ) \big) = 0 \]
for all $s\in(0,1)$, hence we can apply \cite[Proposition 1.62]{AFP00} to infer that
\[ \begin{aligned}
\lim_{r \to 0^+} \frac{ \cH^k \res E (B(p,r) \setminus X(p,\T_p,s))}{\cH^k \res E(B(p,r))} & = \lim_{r \to 0^+} \frac{ (T_{p,r})_\#(\cH^k\res E) (B(e,1) \setminus X(e,\T_p,s))}{\cH^k \res E(B(p,r))} 
\\ & = \cH^k \res \T_p( B(e,1) \setminus X(e,\T_p,s) ) = 0.
\end{aligned} \]
This, together with Lemma \ref{lemma:densities}, shows that at $\cH^k$-a.e. $p \in \G$ the homogeneous subgroup $\T_p$ is a $(k,\G)$-approximate tangent group, which can be also taken to be unique, in view of Proposition~\ref{prop:uniqaptan}.

Proof of ({\romannumeral 4}) $\implies$ ({\romannumeral 1}). 
By compactness of $\cH(\G,k)$, we can find $\V_1,\dots,\V_N \in \cH(\G,k)$ such that for every $\V\in\cH(\G,k)$ there exists $\V_i$ such that $\rho(\V,\V_i)<1/3$. 
Since $\cH^k\res E$ is locally finite, we can define the purely $k$-unrectifiable 
part $E_{\mathrm{pu}}$ of $E$ and consider
\[
C_i= \left\{ p \in E_{\mathrm{pu}} :  \text{there exists }\,
{\rm ap}\Tan_\G^k(E,p)\,\text{and}\, \rho({\rm ap}\Tan_\G^k(E,p),\V_i)<\frac{1}{3}\right\},
\]
that satisfy $E_{\mathrm{pu}} \subseteq Z\cup\bigcup_{i = 1}^N C_i$, 
where $\cH^k(Z)=0$.
So, to prove our claim it suffices that 
\[
\cH^k(C_i) = 0\quad \text{for all $i=1,\ldots,N$.}
\]
We fix one of these integers $i$
and select $p\in C_i$. We wish to prove that 
\begin{equation} \label{eq.5.1}
X\left(p, \V_i^\perp, \epsilon_0\right)\subset \G \setminus X\left(p,\T_p,\epsilon_0 \right),
\end{equation}
up to choosing $\ep_0>0$ suitably small.
We choose 
\[
q \in X(p,\V_i^\perp,\epsilon_0)
\]
and notice that by definition of cone we have
\[
d(p^{-1}q,\V_i^\perp) < \epsilon_0 d(p,q).
\]
It follows from \eqref{eq.lemma.proj3} that $d(p^{-1}q,\V_i^\perp) \geq c_\G \|\pi_{\V_i}(p^{-1}q)\|$ so, ultimately, we get
\beq\label{eq:estpiV}
d(p,q) > \frac{c_\G}{\epsilon_0} \| \pi_{\V_i}(p^{-1}q) \|.
\eeq
To conclude that $q \notin X\left(p,\T_p,\epsilon_0 \right)$,
we need to prove that
\[
d(p^{-1}q,\T_p) \geq \epsilon_0 d(p,q),
\]
but this is not immediate, since $\T_p$ is horizontal and \eqref{eq.lemma.proj3} does not provide directly a good lower estimate of the distance. We start from
the weaker lower estimate \eqref{eq.lemma.proj3}, getting
\[
d(p^{-1}q,\T_p) \geq c_\G \| \pi_{\T_p}(p^{-1}q)^{-1} \pi_{\T_p^\perp}(p^{-1}q) \pi_{\T_p}(p^{-1}q) \|,
\]
so we have to show that
\beq\label{eq:cGestimbelow}
c_\G \| \pi_{\T_p}(p^{-1}q)^{-1} \pi_{\T_p^\perp}(p^{-1}q) \pi_{\T_p}(p^{-1}q) \| \geq \epsilon_0 \|p^{-1}q\|
\eeq
for some suitably small $\epsilon_0 > 0$. Recall that $\T_p$ belongs to $B(\V_i, \frac{1}{3})$, so from definition of distance in the horizontal Grassmannian, we obtain
\[
\|\pi_{\T_p}(p^{-1}q) \| - \| \pi_{\V_i}(p^{-1}q) \| \le 
\| \pi_{\V_i}(p^{-1}q)^{-1} \pi_{\T_p}(p^{-1}q) \| \leq \frac{\|p^{-1}q\|}{3}.
\]
It follows that 
\[ \begin{aligned}
\| \pi_{\T_p}(p^{-1}q)^{-1}\pi_{\T_p^\perp}(p^{-1}q) \pi_{\T_p}(p^{-1}q) \| &
 \geq \|p^{-1}q\| - \| \pi_{\T_p}(p^{-1}q) \|
\\ & \geq d(p,q) - \frac{d(p,q)}{3} - \| \pi_{\V_i}(p^{-1}q) \|
\\ & \geq \left(\frac{2}{3} - \frac{\epsilon_0}{c_\G}\right) d(p,q),
\end{aligned} \]
where the last inequality follows from \eqref{eq:estpiV}.
If we multiply the previous inequalities by $c_\G$, we immediately realize that $\epsilon_0$ can be arbitrarily chosen in the
interval $(0,c_\G/3)$ in order to have 
\eqref{eq:cGestimbelow} established.
This concludes the proof of the inclusion \eqref{eq.5.1}. 

The existence of the $(k,\G)$-approximate tangent group at $p$ yields
\[
\frac{\cH^k \res E( B(p,r) \setminus X(p,\T_p,\ep_0))}{r^k} \to 0\quad\text{as}\quad r\to 0^+
\]
and the point is that $\ep_0$ does not depend on $p$. We fix $\lambda>0$ arbitrarily. 
As a consequence, we obtain
\beq\label{eq:estimFederer}
\cH^k \res E \left(X\left(p,r, \V_i^\perp, \ep_0\right) \right) \leq \lambda r^k \ep_0^k 
\eeq
for all $0<r \leq \delta_{p,\lambda}$ for some 
$\delta_{p,\lambda}> 0$ depending on $p$ and $\lambda$.
For any integer $m\ge1$, we define
\[
C_{im}=\set{q\in C_i: \eqref{eq:estimFederer}\text{ holds for $p=q$ and any $0<r<e^{-m}$}}.
\] 
We observe that $C_i=\bigcup_{m\ge 1}C_{im}$.
The conditions of Theorem~\ref{ppwpod} are satisfied for all points $p\in C_{im}$,
up to choosing a possily smaller $\ep_0<c_\G^3$, therefore
\[
\Theta^{\ast k}(C_{im},\,p) \leq 2(21)^k \lambda.
\]
From \cite[2.10.19(1)]{Federer69}, taking into account the arbitrary choice of $\lambda$,
we get $\cH^k(C_{im})=0$, that immediately leads us to the conclusion of the proof.

\bibliography{References}{}
\bibliographystyle{plain}

\end{document}